\documentclass[12pt]{article}
\usepackage{amsmath,amssymb,amsbsy,amsfonts,amsthm,latexsym,
                        amsopn,amstext,amsxtra,euscript,amscd}

\begin{document}




\newfont{\teneufm}{eufm10}
\newfont{\seveneufm}{eufm7}
\newfont{\fiveeufm}{eufm5}
%
%
\newfam\eufmfam
      \textfont\eufmfam=\teneufm \scriptfont\eufmfam=\seveneufm
      \scriptscriptfont\eufmfam=\fiveeufm
%
%
\def\frak#1{{\fam\eufmfam\relax#1}}
%


\def\bbbr{{\rm I\!R}} 
\def\bbbm{{\rm I\!M}}
\def\bbbn{{\rm I\!N}} 
\def\bbbf{{\rm I\!F}}
\def\bbbh{{\rm I\!H}}
\def\bbbk{{\rm I\!K}}
\def\bbbp{{\rm I\!P}}
\def\bbbone{{\mathchoice {\rm 1\mskip-4mu l} {\rm 1\mskip-4mu l}
{\rm 1\mskip-4.5mu l} {\rm 1\mskip-5mu l}}}
\def\bbbc{{\mathchoice {\setbox0=\hbox{$\displaystyle\rm C$}\hbox{\hbox
to0pt{\kern0.4\wd0\vrule height0.9\ht0\hss}\box0}}
{\setbox0=\hbox{$\textstyle\rm C$}\hbox{\hbox
to0pt{\kern0.4\wd0\vrule height0.9\ht0\hss}\box0}}
{\setbox0=\hbox{$\scriptstyle\rm C$}\hbox{\hbox
to0pt{\kern0.4\wd0\vrule height0.9\ht0\hss}\box0}}
{\setbox0=\hbox{$\scriptscriptstyle\rm C$}\hbox{\hbox
to0pt{\kern0.4\wd0\vrule height0.9\ht0\hss}\box0}}}}
\def\bbbq{{\mathchoice {\setbox0=\hbox{$\displaystyle\rm
Q$}\hbox{\raise
0.15\ht0\hbox to0pt{\kern0.4\wd0\vrule height0.8\ht0\hss}\box0}}
{\setbox0=\hbox{$\textstyle\rm Q$}\hbox{\raise
0.15\ht0\hbox to0pt{\kern0.4\wd0\vrule height0.8\ht0\hss}\box0}}
{\setbox0=\hbox{$\scriptstyle\rm Q$}\hbox{\raise
0.15\ht0\hbox to0pt{\kern0.4\wd0\vrule height0.7\ht0\hss}\box0}}
{\setbox0=\hbox{$\scriptscriptstyle\rm Q$}\hbox{\raise
0.15\ht0\hbox to0pt{\kern0.4\wd0\vrule height0.7\ht0\hss}\box0}}}}
\def\bbbt{{\mathchoice {\setbox0=\hbox{$\displaystyle\rm
T$}\hbox{\hbox to0pt{\kern0.3\wd0\vrule height0.9\ht0\hss}\box0}}
{\setbox0=\hbox{$\textstyle\rm T$}\hbox{\hbox
to0pt{\kern0.3\wd0\vrule height0.9\ht0\hss}\box0}}
{\setbox0=\hbox{$\scriptstyle\rm T$}\hbox{\hbox
to0pt{\kern0.3\wd0\vrule height0.9\ht0\hss}\box0}}
{\setbox0=\hbox{$\scriptscriptstyle\rm T$}\hbox{\hbox
to0pt{\kern0.3\wd0\vrule height0.9\ht0\hss}\box0}}}}
\def\bbbs{{\mathchoice
{\setbox0=\hbox{$\displaystyle     \rm S$}\hbox{\raise0.5\ht0\hbox
to0pt{\kern0.35\wd0\vrule height0.45\ht0\hss}\hbox
to0pt{\kern0.55\wd0\vrule height0.5\ht0\hss}\box0}}
{\setbox0=\hbox{$\textstyle        \rm S$}\hbox{\raise0.5\ht0\hbox
to0pt{\kern0.35\wd0\vrule height0.45\ht0\hss}\hbox
to0pt{\kern0.55\wd0\vrule height0.5\ht0\hss}\box0}}
{\setbox0=\hbox{$\scriptstyle      \rm S$}\hbox{\raise0.5\ht0\hbox
to0pt{\kern0.35\wd0\vrule height0.45\ht0\hss}\raise0.05\ht0\hbox
to0pt{\kern0.5\wd0\vrule height0.45\ht0\hss}\box0}}
{\setbox0=\hbox{$\scriptscriptstyle\rm S$}\hbox{\raise0.5\ht0\hbox
to0pt{\kern0.4\wd0\vrule height0.45\ht0\hss}\raise0.05\ht0\hbox
to0pt{\kern0.55\wd0\vrule height0.45\ht0\hss}\box0}}}}
\def\bbbz{{\mathchoice {\hbox{$\sf\textstyle Z\kern-0.4em Z$}}
{\hbox{$\sf\textstyle Z\kern-0.4em Z$}}
{\hbox{$\sf\scriptstyle Z\kern-0.3em Z$}}
{\hbox{$\sf\scriptscriptstyle Z\kern-0.2em Z$}}}}
\def\ts{\thinspace}

\newtheorem{theorem}{Theorem}
\newtheorem{lemma}[theorem]{Lemma}
\newtheorem{claim}[theorem]{Claim}
\newtheorem{cor}[theorem]{Corollary}
\newtheorem{prop}[theorem]{Proposition}
\newtheorem{definition}{Definition}
\newtheorem{question}[theorem]{Open Question}

\def\squareforqed{\hbox{\rlap{$\sqcap$}$\sqcup$}}
\def\qed{\ifmmode\squareforqed\else{\unskip\nobreak\hfil
\penalty50\hskip1em\null\nobreak\hfil\squareforqed
\parfillskip=0pt\finalhyphendemerits=0\endgraf}\fi}




\newcommand{\ignore}[1]{}

\newcommand{\comm}[1]{\marginpar {\fbox{#1}}}

\hyphenation{re-pub-lished}

\baselineskip 15pt

\def\lln{{\mathrm Lnln}}
\def\ad{{\mathrm ad}}

\def \F{{\bbbf}}
\def \K{{\bbbk}}
\def \Z{{\bbbz}}
\def \N{{\bbbn}}
\def \Q{{\bbbq}}
\def \R{{\bbbr}}
\def\Fp{\F_p}
\def \fp{\Fp^*}

\def\cA{{\mathcal A}}
\def\cE{{\mathcal E}}
\def\cI{{\mathcal I}}
\def\cQ{{\mathcal Q}}
\def\cR{{\mathcal R}}
\def\cX{{\mathcal X}}
\def\cU{{\mathcal U}}
\def\Zm{\Z_m}
\def\Zt{\Z_t}
\def\Zp{\Z_p}
\def\Fp{\F_p}
\def\Um{\Uc_m}
\def\Ut{\Uc_t}
\def\Up{\Uc_p}

\def\Ikl{\widetilde{I}_\ell(K,H,N)}
\def\\{\cr}
\def\({\left(}
\def\){\right)}
\def\fl#1{\left\lfloor#1\right\rfloor}
\def\rf#1{\left\lceil#1\right\rceil}
\def\ep{\mathbf{e}}

\newcommand{\vt}{\vartheta}
\newcommand{\floor}[1]{\lfloor {#1} \rfloor }


\title{Character Sums and Congruences with $n!$}

\author{
{Moubariz~Z.~Garaev}\\
\normalsize{Instituto de Matem{\'a}ticas,  Universidad Nacional Aut\'onoma de
M{\'e}xico}
\\
\normalsize{C.P. 58180, Morelia, Michoac{\'a}n, M{\'e}xico} \\
\normalsize{\tt garaev@matmor.unam.mx} \\
\and
{ Florian~Luca} \\
\normalsize {Instituto de Matem{\'a}ticas, Universidad Nacional Aut\'onoma de
M{\'e}xico}
\\
\normalsize{C.P. 58180, Morelia, Michoac{\'a}n, M{\'e}xico} \\
\normalsize{\tt fluca@matmor.unam.mx} \\
\and
{Igor E.~Shparlinski} \\
\normalsize{Department of Computing, Macquarie University} \\
\normalsize{Sydney, NSW 2109, Australia} \\
\normalsize{\tt igor@ics.mq.edu.au}}
\date{}

\pagenumbering{arabic}

\maketitle


\begin{abstract} We estimate character sums with $n!$, on average, and
individually. These bounds are used to derive new results about various
congruences modulo a prime $p$ and obtain  new information about
the spacings between quadratic nonresidues  modulo $p$. 
In particular, we show that there exists
a positive integer $n\ll p^{1/2+\varepsilon}$, such
that $n!$ is a primitive
root modulo $p$. We also show that every nonzero congruence class
$a \not \equiv 0 \pmod p$
can be represented as a product of 7 factorials,
$a \equiv n_1! \ldots n_7! \pmod p$, where $\max \{n_i \ |\ i=1,\ldots
7\}=O(p^{11/12+\varepsilon})$, and we find the asymptotic formula
for the number of such representations.  
Finally, we show that products of 4 factorials $n_1!n_2!n_3!n_4!,$ with 
$\max\{n_1, n_2, n_3, n_4\}=O(p^{6/7+\varepsilon})$ represent ``almost all''
residue classes
modulo p, and that products of 3 factorials $n_1!n_2!n_3!$ with $\max\{n_1,
n_2, n_3\}=O(p^{5/6+\varepsilon})$ are uniformly  distributed modulo $p$.
\end{abstract}


\section{Introduction}

Throughout this paper, $p$ is an odd prime. Very little
seems to be known about the distribution of $n!$ modulo $p$.
In {\bf F11} in~\cite{RKG}, it is conjectured that
about $p/e$ of the residue classes $a\pmod p$ are missed
by the sequence $n!$. If this were so, the sequence
$n!$ modulo $p$ should assume about $(1-1/e)p$ distinct values. Some
results in this spirit appear in~\cite{CVZ}. The above conjecture
    immediately implies  that every residue class
$a$ modulo $p$ can be represented as a product of at most two factorials.
Unconditionally, it is easy to see that three factorials suffice. Indeed,
$0\equiv p!\pmod p$, and, as it has been remarked in~\cite{ErSt}, equation~(5),
see also~\cite{CVZ},
the {\it Wilson theorem\/} implies that
\begin{equation}
\label{Wilson Thm}
b!\cdot (p-1-b)!\equiv  (-1)^{b+1} \pmod p
\end{equation}
holds for any $b \in \{1,\ldots,p-1\}$. Therefore, if $a\in \{1,\ldots,p-1\}$,
then with $b \equiv a^{-1} \pmod p$, we have
\begin{equation}
\label{Wilson Repr}
a \equiv ((p-1)!)^{r_{b}}(b-1)!\cdot (p-1-b)!  \pmod p,
\end{equation}
where
$r_b\in \{0,1\}$ is such that $r_b \equiv b+1\pmod 2$.
However, the above argument does not apply to proving the
existence of representations involving factorials of integers of
restricted size,
neither can it be used for estimation of the number of representations.

In this paper, we first estimate character sums with $n!$ on the average,
and individually. We use these estimates to show that for every
$\varepsilon$ and $p$ sufficiently large, there exists a value of $n$
with $n=O(p^{1/2+\varepsilon})$ and such that $n!$ is a primitive root modulo
$p$.

We  apply these estimates to prove  that every residue class
$a\not\equiv 0 \pmod
p$,  can be represented as a product of $7$ factorials, $a\equiv n_1!\ldots
n_7! \pmod p$ with $\max \{n_i \ | \ i=1,\ldots, 7\}\ll p^{11/12+\varepsilon}
$. If we only want that ``most'' of the residue classes modulo $p$ be
represented as a product of factorials in the same range as above
(and even a slightly better one), then we show that four
factorials suffice.
Moreover, our results imply that for every $\varepsilon>0$ and sufficiently
large $p$,
every residue class $a\not\equiv 0 \pmod p$ can be represented as a product of
$\ell = \fl{\varepsilon^{-1}} + 5$  factorials,
$a\equiv n_1!\ldots n_\ell ! \pmod p$,  where
$\max \{n_i \ | \ i=1,\ldots,\ell\}\ll p^{1/2+\varepsilon}$.
We also show that products
of three factorials $n_1! n_2! n_3!$, with
$\max \{n_1, n_2, n_3\}  = O( p^{5/6+\varepsilon})$, are
uniformly distributed modulo $p$.

Our basic tools are the {\it Weil bound\/} for character sums,
see~\cite{Li,LN,W},
   and the {\it Lagrange theorem\/} bounding the number of zeros of
a non-zero polynomial over a field.

Some of the results of this paper have 
found 
applications to the study of arithmetic
properties of expressions of the form $n!+f(n)$, where $f(n)$ is 
a polynomial with integer coefficients (see ~\cite{LSh1}), 
or a linearly recurrent sequence of integers (see ~\cite{LSh2}).
In particular, an improvement of a result of Erd{\H o}s  and  Stewart~\cite{ErSt},
obtained in~\cite{LSh1},  is based on these results.

Throughout the paper the  implied constants in symbols `$O$' and  `$\ll$' 
may occasionally, where obvious, depend on   integer parameters
$\ell$  and $d$ and  a small real parameter $\varepsilon > 0$,
 and are absolute otherwise
(we recall that $A \ll B$  is equivalent
to $A = O(B)$).

\bigskip

{\bf Acknowledgements.} The authors would like to thank
Vsevolod Lev  for several useful comments.
During the preparation of this paper,
F.~L.\ was supported in part by grants
SEP-CONACYT 37259-E and 37260-E, and
I.~S.\ was supported in part by ARC grant  DP0211459.

\section{Character Sums}

Let $\F_p$  be a finite field  of
$p$ elements. We always assume that $\F_p$ is represented by the elements
of the set $\{0, 1, \ldots, p-1\}.$

Let $\cX$ denote the set of  multiplicative characters of the multiplicative
group $\F_p^*$
     and let $\cX^{*} = \cX
\backslash \{\chi_0\}$  be   the set   of nonprincipal characters.

We also define
$$
\ep(z) = \exp(2 \pi i z/p),
$$
which is an additive character  of $\F_p$.

It is useful to recall the identities
$$
\sum_{\chi\in \cX} \chi( u) =
\left\{ \begin{array}{ll}
0,& \quad \mbox{if}\ u\not \equiv 1 \pmod p, \\
p-1,& \quad \mbox{if}\ u \equiv 1 \pmod p,
\end{array} \right.
$$
and
$$
\sum_{a=0}^{p-1} \ep( au) =
\left\{ \begin{array}{ll}
0,& \quad \mbox{if}\ u\not \equiv 0 \pmod p, \\
p,& \quad \mbox{if}\ u \equiv 0 \pmod p,
\end{array} \right.
$$
which we will repeatedly use, in particular to relate the  number of solutions
of various congruences and character sums.

Given $\chi \in \cX$,  a polynomial
$f \in \F_p[X]$,   and an element $a \in \F_p$,
we consider character sums
$$
T(\chi,f,H,N)= \sum_{n = H+1}^{H+N} \chi\(n!\)
\ep(f(n))
$$
where we simply write $T(\chi,H,N)$ is $f$ is identical to zero, and
$$
S(a,H,N)= \sum_{n = H+1}^{H+N} \ep\(a n!\) .
$$

We obtain a nontrivial upper bound for ``individual'' sums $T(\chi,f,H,N)$,
and also nontrivial upper bounds for the moments of
$T(\chi,f,H,N)$ and $S(a,H,N)$.

\begin{theorem}
\label{thm:Char-Indiv}
Let $H$ and $N$ be integers with  $0 \le H< H+ N < p$.
Then for any fixed   integer $d \ge 0$,
the following bound holds:
$$
\max_{\deg f = d} \, \max_{\chi \in \cX^*} \, |T(\chi,f,H,N)| \ll
N^{3/4} p^{1/8} (\log p)^{1/4}.
$$
\end{theorem}

\begin{proof} For any integer $k \ge 0$ we have
$$
T(\chi,f,H,N) = \sum_{n =H+1}^{H+N}
\chi\((n+k)!\)  \ep(f(n+k))+ O(k).
$$
Therefore, for any integer $K \ge 0$,
\begin{equation}
\label{eq:T and W}
T(\chi,f,H,N) = \frac{1}{K} W + O(K),
\end{equation}
where
\begin{eqnarray*}
W & = & \sum_{k = 0}^{K-1}\sum_{n =H+1}^{H+N}
     \chi\((n+k)!\)  \ep(f(n+k))\\ &= & \sum_{n =H+1}^{H+N}
     \sum_{k = 1}^{K} \chi\( n!
\prod_{i=1}^k (n+ i)\) \ep(f(n+k))\\ & = & \sum_{n =H+1}^{H+N}
\chi\( n!\) \sum_{k = 0}^{K-1}
\chi\( \prod_{i=1}^k (n+ i)\) \ep(f(n+k)).
\end{eqnarray*}
Applying the Cauchy  inequality, we derive
\begin{eqnarray*}
|W|^2 &\le & N  \sum_{n =H+1}^{H+N}  \left| \sum_{k = 0}^{K-1}
\chi\(  \prod_{i=1}^k (n + i)\) \ep(f(n+k))\right|^{2} \\
     &  =  & N \sum_{k, m= 0}^{K-1}\sum_{n
=H+1}^{H+N}\hskip-20pt{\phantom{\sum}}^* 
\chi\(  \Psi_{k,m}(n)\) \ep(f(n+k)-f(n+m)),
\end{eqnarray*}
where
$$
\Psi_{k,m}(X)  =
\prod_{i=1}^{k} (X + i)  \prod_{j=1}^{m} (X + j)^{-1}
$$
and hereafter $\sum{}^* $ means that the poles of the corresponding rational
function are excluded from the summation (we also recall
that  
$|z|^2 = z
\overline{z}$ for any complex number
$z$, and that $  \overline{\chi}(a) = \chi(a^{-1})$  holds for every integer
$a \not
\equiv 0 \pmod p$ where $\overline{\chi}$ is
the conjugate character  of ${\chi}$).

Clearly, if $K < p$ then,  unless $k = m$, the rational function
$\Psi_{k,m}(X)  $,  has at least one
single root or pole, and thus  is not a power of any other rational
function modulo $p$.

For the  $O(K)$ choices of   $0 \le k =m\le K-1$,
     we estimate the sum over $n$ trivially as $N$.

For the other  $O(K^2)$  choices of  $0 \le k,m\le K-1$,
using  the Weil bound,  given in
Example~12 of Appendix~5 of~\cite{W}
(see also Theorem~3 of Chapter~6 in~\cite{Li}, or
Theorem~5.41 and the comments to Chapter~5 of~\cite{LN}),
we see that, because $\chi \in \cX^*$, then, for any $a \in \F_p$, we have
$$
\sum_{n =0}^{p-1}\hskip-20pt{\phantom{\sum}}^*   \chi\(  \Psi_{k,m}(n)\)
\ep(f(n+k)-f(n+m) + an)
      = O(Kp^{1/2 }).
$$
Therefore, by the standard reduction of incomplete sums to complete ones,
(see~\cite{Chalk}),
we deduce
$$
\sum_{n =H+1}^{H+N} \hskip-20pt{\phantom{\sum}}^*  \chi\(  \Psi_{k,m}(n)\)
\ep(f(n+k)-f(n+m))
      = O(Kp^{1/2}\log p).
$$
Putting everything together, we get
$$
W^2 \ll  N \(K N +  K^{3} p^{1/2} \log p\).
$$
Therefore, by~\eqref{eq:T and W}, we derive
$$
T(\chi,f,H,N) \ll   N  K^{-1/2} +
     K^{1/2} N^{1/2} p^{1/4} (\log p)^{1/2}
+  K.
$$
Taking $K = \fl{ N^{1/2} p^{-1/4} (\log p)^{-1/2} }$, we  finish the proof.
\end{proof}

It is clear that for any $\varepsilon > 0$ there exists some $\delta >0$,
such that if $N \ge p^{1/2 + \varepsilon}$   then
$$
|T(\chi,f,H,N)| \le  N p^{-\delta}.
$$
provided that $p$ is large enough.

Clearly, Theorem~\ref{thm:Char-Indiv}  immediately implies that 
among the values of $n!$,  where 
$n =H+1, \ldots, H+N$, 
there are  $N/2 + O\(N^{3/4} p^{1/8} (\log p)^{1/4}\)$ quadratic 
residues and nonresidues.   Remarking that each change in the value 
of the Legendre symbol $(n!/p)$ corresponds to a quadratic non-residue
$n$ we can derive a certain result about the distribution of spacings between
quadratic non-residues modulo $n$ which does not seem to follow from 
any of the previously known results, see~\cite{KurlRud}.

Let $n_j$ be the $j$th quadratic nonresidue modulo $p$ and
let  $d_j = n_{j} - n_{j-1}$, the   $j$th spacing,  
$j =1, \ldots, (p-1)/2$, where we put $n_0 = 0$.

\begin{cor} Let $J$ be an integer with $p^{1/2} \log p\le J \le (p-1)/2$.
Then the following bound holds:
$$
\sum_{j=1}^J (-1)^j d_j   \ll  J^{3/4} p^{1/8}(\log p)^{1/4}.
$$
\end{cor}

\begin{proof}
We have 
$$
\sum_{j=1}^J (-1)^{j-1} d_j  =  \sum_{n = 0}^{n_J-1} \(\frac{n!}{p}\)
$$
 {From} the {\it  Polya--Vinogradov bound\/}
\begin{equation}
\label{eq:Vin-Pol}
\max_{\chi \in \cX^*} \max_{0 \le h \le k \le p-1}
\left| \sum_{c = h+1}^k  \chi(c) \right|  \ll p^{1/2} \log p
\end{equation}
we see that  $n_J= 2 J +  O(p^{1/2} \log p) \ll J$
 and by Theorem~\ref{thm:Char-Indiv}
we derive the result.
\end{proof}

Obviously 
$$
\sum_{j=1}^J  d_j  =  n_J = 2 J +  O(p^{1/2} \log p),
$$
which demonstrates that for every $J \ge p^{1/2 + \varepsilon}$
the odd and even spacings $d_j$, $j=1, \ldots, J$, are 
of approximately the same total length.

We now denote by $Q(H,N)$  the number of  $n =H+1, \ldots, H+N$
such that $n!$  is a primitive root modulo $p$.

\begin{cor} Let $H$ and $N$ be integers with $0 \le H< H+ N < p$.
Then, for any fixed   $\varepsilon > 0$,
the following bound holds:
$$
Q(H,N) = N  \frac{\varphi(p-1)}{p-1} 
+O\( N^{3/4} p^{1/8+ \varepsilon } \).
$$
\end{cor}

\begin{theorem}
\label{thm:Char-Aver}
Let $H$ and $N$ be integers with  $0 \le H< H+ N < p$.
Then for any fixed  integer $d \ge 0$,
the following bound holds:
$$
\max_{\deg f = d} \, \sum_{\chi \in \cX} \, |T(\chi,f,H,N)|^{2} \ll
p N^{3/2}  .
$$
\end{theorem}

\begin{proof} Arguing as in the proof of Theorem~\ref{thm:Char-Indiv},
and applying the
H{\"o}lder  inequality  to~\eqref{eq:T and W}, we derive that for any $K$
\begin{eqnarray*}
|T(\chi,f,H,N) |^{2} & \ll  &K^{-2} N \sum_{k, m= 0}^{K-1}
\sum_{n =H+1}^{H+N}  \hskip-20pt{\phantom{\sum}}^*  \chi\(  \Psi_{k,m}(n)\) \\
& & \qquad  \qquad  \qquad  \qquad  \quad \ep\( f(n+k)-f(n+m)\)  + K^{2},
\end{eqnarray*}
where
$$
\Psi_{km}(X)  =\prod_{i=1}^{k} (X  + i)
\prod_{j=1}^{m} (X + j)^{-1}
$$
and  as before $\sum{}^* $ means that the poles of the corresponding rational
function are excluded from the summation.   Therefore,
\begin{eqnarray*}
\lefteqn{\sum_{\chi \in \cX} \, |T(\chi,f,H,N)|^{2} }\\
& &\qquad  \ll
K^{-2} N\sum_{k,m= 0}^{K-1} \sum_{n =H+1}^{H+N} 
\hskip-20pt{\phantom{\sum}}^*    \sum_{\chi \in \cX}
\chi\(  \Psi_{k,m}(n)\) + pK^{2}.
\end{eqnarray*}
The sum over $\chi$ vanishes, unless
\begin{equation}
\label{eq:Psi=1}
\Psi_{k,m}(n) \equiv 1 \pmod p
\end{equation}
in which case it is equal to $p-1$.
For $K$ pair $(k,m)$ with $k=m$ then there are $N$ possible solutions
to~\eqref{eq:Psi=1},
for other $O(K^2)$ pairs there are at most $K$ solutions to~\eqref{eq:Psi=1}. 
 Thus
\begin{eqnarray*}
\sum_{\chi \in \cX} \, |T(\chi,f,H,N)|^{2}  &\ll  &
K^{-2} N \(  K^{3} + K  N\)p + pK^{2}\\
& = &\(N K  + N^{2}K^{-1} +  K^{2}\) p.
\end{eqnarray*}
Taking $K = \fl{ N^{1/2}}$,  we  finish the proof.
\end{proof}

\begin{theorem}
\label{thm:Exp-Aver}
Let $H$ and $N$ be integers with  $0 \le H< H+ N < p$.
Then for any fixed   integer $\ell \ge 1$,
the following bound holds:
$$
\sum_{a =0}^{p-1} |S(a,H,N)|^{2} \ll
p N^{3/2}  .
$$
\end{theorem}

\begin{proof} Arguing as in the proof of Theorem~\ref{thm:Char-Indiv},
we derive
$$\sum_{a =0}^{p-1}|S(a,H,N)|^{2}
\ll  K^{-2}
N\sum_{k,m= 0}^{K-1}
\sum_{n =H+1}^{H+N}  \sum_{a =0}^{p-1} \ep\( a n! \Phi_{k,m}(n)\) + p
K^{2},
$$
where
$$
\Phi_{k,m}(X)  =
\prod_{i=1}^{k_\nu} (X  + i)  -
\prod_{j=1}^{k_{\nu+s}} (X+ j).
$$
The sum over $a$ vanishes, unless
\begin{equation}
\label{eq:Phi=1}
n!\Phi_{k,m}(n) \equiv  0 \pmod p,
\end{equation}
in which case it is equal to $p$.

As before,
we see that  $\Phi_{k,m}(X)$ is a nonconstant polynomial of degree $O(K)$,
unless $k=m$
 Because $n! \not\equiv 0
\pmod p$ for
$0 \le H < n \le H + N < p$, we derive
$$
\sum_{a =0}^{p-1}|S(a,H,N)|^{2} \ll   \(N K  +
N^{2}K^{-1} +
K^{2}\) p.
$$
Taking $K = \fl{ N^{1/2}}$ and remarking that with this value
of $K$ the last term never dominates,  we  finish the proof.
\end{proof}

\section{Sums and Products of Factorials}

For  integer $\ell\ge 1 $ and  $H$ and $N$  with  $0 \le H< H+ N < p$
let us denote by $I_\ell(H,N)$ and  $J_\ell(H,N)$ the number of solutions
to the congruences
$$
\prod_{i=1}^\ell  n_i! \equiv \prod_{i=\ell+1}^{2\ell} n_i!  \pmod p,
\qquad  H +1 \le n_1, \ldots, n_{2\ell}\le H+N,
$$
and
$$
\sum_{i=1}^\ell  n_i! \equiv \sum_{i=\ell+1}^{2\ell} n_i!  \pmod p,
\qquad  H +1 \le n_1, \ldots, n_{2\ell} \le H+N,
$$
respectively. 

{From} the properties of multiplicative and additive characters we immediately
conclude that 
 \begin{equation}
\label{eq:T and I}
\frac{1}{p-1}\sum_{\chi \in \cX}  |T(\chi,f,H,N)|^{2\ell} \le  \frac{1}{p-1}
\sum_{\chi
\in
\cX} |T(\chi,H,N)|^{2\ell} = I_\ell(H,N)  
\end{equation}
and 
 \begin{equation}
\label{eq:S and J}
 \frac{1}{p} \sum_{a=0}^{p-1} \left| S(a,H,N)\right|^{2\ell} = J_\ell(H,N) .
\end{equation}

The same arguments as in the proof of Theorem~\ref{thm:Exp-Aver}
lead to the bound 
$$
 J_\ell(H,N) \ll N^{2\ell - 1 + 1/(\ell + 1) } .
$$
We now show that for $I_\ell(H,N) $ one can derive a more precise estimate.

\begin{theorem}
\label{thm:AverValues}
Let $H$ and $N$ be integers with  $0 \le H< H+ N < p$.
Then for any fixed   integer  $\ell \ge 1$,
the following bound holds:
$$
I_\ell(H,N) \ll N^{2\ell - 1 + 2^{-\ell}}  .
$$
\end{theorem}

\begin{proof} We prove this bound by induction.
If $\ell =1$ then Theorem~\ref{thm:Char-Aver}  taken with $f(X) = 0$,
together with~\eqref{eq:T and I} 
immediately imply  the desired bound
$I_1(H,N) \ll N^{3/2}$.

Now assume that $\ell \ge 2$  and that 
 $I_{\ell-1}(H,N) \ll N^{2\ell - 3 + 2^{-\ell+1}}$.  
We fix some
$K< N$ and note that by the Cauchy inequality we have
\begin{eqnarray*}
\left| \sum_{n=H+1}^{H+N}
\chi(n!)\right|^{2} &= &\left| \sum_{k=1}^K  \sum_{H+ (k-1)N/K < m \le H+kN/K}
\chi(n!)\right|^{2}\\
& \le & 
K \sum_{k=1}^K  \left|\sum_{H+ (k-1)N/K < m \le H+kN/K}
\chi(n!)\right|^{2} .
\end{eqnarray*}
Therefore
\begin{eqnarray*}
I_\ell(H,N) & = & \frac{K}{p-1} \sum_{k=1}^K   \sum_{\chi \in \cX} 
\left|\sum_{H+ (k-1)N/K < m \le H+kN/K}
\chi(n!)\right|^{2}  \left| \sum_{n=H+1}^{H+N}
\chi(n!)\right|^{2\ell-2}\\
& = & K  \Ikl,
\end{eqnarray*}
where  $\Ikl$ is the number of solutions to the congruence
$$
m_1! \prod_{i=1}^{\ell-1}  n_i! \equiv m_2! \prod_{i=\ell}^{2\ell-2} n_i!  \pmod p
$$
with $H +1 \le n_1, \ldots, n_{2\ell-2}\le H+N$  and 
$H+ (k-1)N/K < m_1,m_2 \le H+kN/K$ for some $k =1 , \ldots, K$.
For each of $N$ pairs $(m_1, m_2)$ with 
 $m_1 = m_2$, there are  exactly $I_{\ell-1}(H,N)$ 
solutions. Also we see that if  $n_1, \ldots, n_{2\ell-2}$
are given  then for each fixed value of $r=m_1-m_2$,
there are no more then $|r|$ values solutions in $m_1, m_2$
(because at least one of $m_1, m_2$ satisfies a 
nontrivial polynomial  congruence of degree $|r|$). Certainly $r = O(N/H)$.
Putting everything together and using the induction
assumption we obtain
$$
\Ikl \ll N I_{\ell-1}(H,N) + (N/K)^2 N^{2\ell -2} = 
N^{2\ell -2 + 2^{-\ell +1}}  +  N^{2\ell} K^{-2}.
$$
Therefore $I_\ell(H,N)  \ll K N^{2\ell -2 + 2^{-\ell +1}}  +  N^{2\ell} K^{-1}$.
Choosing  $K = \rf{ N^{1-2^{-\ell}}}$,
 we obtain the desired bound.  \end{proof}

We now show that, for $N \ge p^{1/2 + \varepsilon}$ the above  bound
on $I_\ell(H,N)$,
combined with Theorem~\ref{thm:Char-Indiv}, produces an asymptotic formula
for  $I_\ell(H,N)$.  In particular for $H=0$, $N = p-1$, this
asymptotic formula
is nontrivial for $\ell \ge 4$.

\begin{theorem}
\label{thm:Two Products}
Let $H$ and $N$ be integers with  $0 \le H< H+ N < p$.
Then for any fixed   integers  $\ell \ge r\ge 1$,
the following bound holds:
$$
I_\ell(H,N)  = \frac{N^{2\ell}}{p-1} + O\(N^{3\ell/2+r/2 -1 +2^{-r}}
p^{(\ell-r)/4}
(\log p)^{(\ell-r)/2 }\).
$$
\end{theorem}

\begin{proof} Similar to~\cite{Kar}, we have 
 \begin{eqnarray*}
\lefteqn{I_\ell(H,N) =\frac{1}{p-1} \sum_{\chi \in \cX} \left| \sum_{n=H+1}^{H+N}
\chi(n!)\right|^{2\ell} }\\
& & \qquad=   \frac{N^{2\ell}}{p-1} +
\frac{1}{p-1}
\sum_{\chi \in \cX^*} \left| \sum_{n=H+1}^{H+N} \chi(n!)\right|^{2\ell} \\
    &  & \qquad = \frac{N^{2\ell}}{p-1} + O\( \frac{1}{p-1}
\max_{\chi \in \cX^*} \left| \sum_{n=H+1}^{H+N}
\chi(n!)\right|^{2\ell - 2r}
\sum_{\chi \in \cX} \left| \sum_{n=H+1}^{H+N} \chi(n!)\right|^{2r}\),
\end{eqnarray*}
(note that  in the last sum we bring back the term corresponding to $\chi=
\chi_0$).
The result follows from~\eqref{eq:T and I} and 
Theorems~\ref{thm:Char-Indiv} and~\ref{thm:AverValues}.
\end{proof}

In particular, using Theorem~\ref{thm:Two Products}
 with $r=1$ we obtain
$$
I_4(0,p-1)  = p^7 \( 1 + O(p^{-1/4} (\log p)^{3/2})\).
$$

We now  denote by $F_\ell(a,H,N)$  the number of solutions
to the congruence
$$
\prod_{i=1}^\ell  n_i! \equiv a \pmod p, \qquad  H +1 \le n_1, \ldots, n_{\ell}
\le H+N,
$$
where $a \in \F_p^*$.

The same arguments as the ones used
in the proof of Theorem~\ref{thm:Two Products}
imply:

\begin{theorem}
\label{thm: One Product}
Let $H$ and $N$ be integers with  $0 \le H< H+ N < p$.
Then, for any fixed integers  $\ell \ge 2r\ge 1$,
the following bound holds:
$$
\max_{a \in \F_p^*} \left|F_\ell(a,H,N) - \frac{N^{\ell}}{p-1}\right|
\ll  N^{3\ell/4 +  r/2  - 1  +2^{-r}} p^{(\ell - 2r)/8} (\log p)^{(\ell - 2r)/4}.
$$
\end{theorem}

In particular, using Theorem~\ref{thm: One Product} with $r=1$
we obtain
$$
F_7(a,0,p-1) = p^6( 1 + O\(  p^{-1/8} (\log p)^{5/4}\),
$$
and for any $\varepsilon > 0$, using Theorem~\ref{thm: One Product} with $r=2$
we obtain
$$
F_7(a,H,N) = \frac{N^{7}}{p} ( 1 +o(1)), \qquad \text{for} \
N \ge   p^{11/12+\varepsilon},
$$
hold for all $a \in \F_p^*$.

Let $V_\ell(H,N)$ be the number of $a \in \F_p^*$  for which
$F_\ell(a,H,N) >0$, that is,
$$
V_\ell(H,N) = \# \left\{ \prod_{i=1}^\ell  n_i! \pmod p, \ | \  H +1
\le n_1, \ldots,
n_{\ell}
\le H+N\right\}.
$$

\begin{theorem}
\label{thm:ValueSet-1}
Let $H$ and $N$ be integers with  $0 \le H< H+ N < p$.
Then for any fixed   integers  $\ell \ge r\ge 1$,
the following bound holds:
$$
V_\ell(H,N) = p  + O\( N^{-\ell/2  +r/2  -1 +2^{-r} } p^{(\ell- r + 8)/4}
(\log p)^{(\ell-r)/2}
\).
$$
\end{theorem}

\begin{proof} We may assume that $\ell\ge 2$, otherwise there is nothing 
to prove. Let
$$
\cE=\left\{h  \in \F_p \ | \    h \not \equiv  \prod_{i=1}^\ell  n_i!
\pmod{p},
\ H +1
\le n_1, \ldots, n_{\ell}
\le H+N\right\}.
$$
Then,
$$
\frac{1}{p-1}\sum_{\chi \in \cX}\sum_{n_1, \ldots, 
n_\ell=H+1}^{H+N}\sum_{h\in \cE}
\chi(n_1!\ldots n_{\ell}!h^{-1})=0.
$$
Separating  the term corresponding to $\chi_0$ and, for $\chi\in 
\cX^*$, applying
Theorem~\ref{thm:Char-Indiv}  to the sums over
$n_1, \ldots, n_{\ell-r}$, we obtain
$$
\frac{\# \cE N^{\ell}}{p-1}\le \left(N^{3/4}p^{1/8}\log^{1/4}
p\right)^{\ell -r}
\frac{1}{p-1}\sum_{\chi \in \cX^*}\left|\sum_{h\in \cE} \chi(h)\right|  \left|
\sum_{n 
=H+1}^{H+N}\chi(n! )\right|^r.
$$
As before, 
we now extend summation over all characters $\chi \in \cX$ and by the  Cauchy
inequality, we  derive from~\eqref{eq:T and I} and Theorem~\ref{thm:AverValues}
\begin{eqnarray*}
 \lefteqn{\(\sum_{\chi \in \cX}\left| \sum_{h\in \cE} \chi(h)\right|  \left|
\sum_{n
=H+1}^{H+N}\chi(n! )
\right|^r\)^2   \le  \sum_{\chi \in \cX} \left| \sum_{h\in 
\cE}\chi(h^{-1}) \right|^2
\sum_{\chi \in \cX} \left|\sum_{n =H+1}^{H+N}\chi(n!) \right|^{2r} }\\
& & \qquad \qquad \qquad \qquad \qquad
  = (p-1) I_r(H,N) \# \cE   \ll  (p-1)^2  N^{2r -1 + 2^{-r}} \# \cE  .
\end{eqnarray*}
Therefore,
$$
\frac{ N^{\ell}\# \cE}{p-1}\le \left(N^{3/4}p^{1/8}\log^{1/4}
p\right)^{\ell-r}\cdot \left
( \# \cE N^{2r -1 + 2^{-r}}\right)^{1/2},
$$
which finishes the proof.
\end{proof}

In particular, using Theorem~\ref{thm:ValueSet-1}
with $r=2$ we see that for $N>p^{6/7+\varepsilon}$,   we have
   that only $o(p)$
residue classes modulo $p$ cannot be represented as
$n_1!n_2!n_3!n_4!\pmod{p}$ with $H+ 1\le n_1, n_2, n_3,
n_4\le H+N$.

We recall that the {\it discrepancy\/} $D$ of a sequence of
$M$ points $(\gamma_j)_{j=1}^M$
of the unit interval $[0,1]$ is  defined   as
$$
D =\sup_{\cI} \left|\frac{A(\cI) } {M}- |\cI|\right|,
$$
where the supremum is taken over the interval
$\cI = [\alpha, \beta]  \subseteq [0,1]$ of length
$|\cI| = \beta - \alpha$
and $A(\cI)$ is the number of points of this set which belong to $\cI$
(see~\cite{DrTi,KuipNied}).

For  an integer $a$ with $\gcd(a,p) = 1$,  we   denote by
$D_\ell(a,H,N)$  the
discrepancy  of the sequence of fractional parts
$$
\left\{\frac{a}{p}  \prod_{i=1}^\ell  n_i!\right\}, \qquad  H +1 \le
n_1, \ldots, n_{\ell}
\le H+N.
$$

Obviously,
\begin{equation}
\label{eq:D and F}
D_\ell(a,H,N) =
\max_{0 \le h \le k \le p-1} \left| \frac{1}{N^\ell}\sum_{c= h+1}^k
F_\ell(a^{-1} c,H,N) -
\frac{k-h}{p}
\right| + O(p^{-1}),
\end{equation}
thus Theorem~\ref{thm: One Product} can be used to estimate $D_\ell(a,H,N)$.
However, we show that the  Polya--Vinogradov bound~\eqref{eq:Vin-Pol}
leads to stronger results.

\begin{theorem}
\label{thm:Discrep}
Let $H$ and $N$ be integers with  $0 \le H< H+ N < p$.
Then for any fixed   integers  $\ell \ge 2r \ge 1$,
the following bound holds:
$$
\max_{1 \le a \le p-1} |D_\ell(a,H,N)|  \ll N^{-\ell/4+r/2 - 1 +2^{-r}} 
p^{(\ell -2r+4)/8} (\log p)^{(\ell -2r+ 4)/4}.
$$
\end{theorem}

\begin{proof} We have
\begin{eqnarray*}
\lefteqn{
\sum_{c = h+1}^k  F_\ell(a^{-1} c,H,N) -
\frac{k-h}{p} N^\ell  }\\
&  & \qquad \qquad =
\frac{1}{(p-1)} \sum_{\chi \in \cX^*}  \sum_{c = h+1}^k
\sum_{n_1, \ldots, n_\ell=H+1}^{H+N}
\chi\(a c^{-1} \prod_{i =1}^\ell n_i!\)\\
&  & \qquad \qquad =
\frac{1}{(p-1)} \sum_{\chi \in \cX^*} \chi(a) \sum_{a = h+1}^k
\overline{\chi}(c)
\(\sum_{n=H+1}^{H+N}
\chi\( n!\)\)^\ell.
\end{eqnarray*}
Thus, applying the bound~\eqref{eq:Vin-Pol}, we  deduce
\begin{eqnarray*} 
\lefteqn{
\left|\sum_{c = h+1}^k  F_\ell(a^{-1} c,H,N) -
\frac{k-h}{p} N^\ell \right|}\\
&  &  \qquad \ll
\max_{\chi \in \cX^*} \left| \sum_{n=H+1}^{H+N}
\chi(n!)\right|^{\ell - 2r}
\sum_{\chi \in \cX} \left| \sum_{n=H+1}^{H+N} \chi(n!)\right|^{2r}
p^{-1/2} \log p,
\end{eqnarray*}
and the result follows from~\eqref{eq:T and I} and Theorems~\ref{thm:Char-Indiv}
and~\ref{thm:AverValues}.
\end{proof}

In particular, using Theorem~\ref{thm:Discrep} with $r=1$ we obtain
that
$$
\max_{1 \le a \le p-1} |D_3(a, 0, p-1)| = O\(p^{-1/8} (\log p)^{5/4}\),
$$
and also that for any $\varepsilon > 0$, 
$$
\max_{1 \le a \le p-1} |D_3(0, p-1)| = o(1), \qquad \text{for} \
N \ge p^{5/6+\varepsilon}.
$$

We also note that Theorem~\ref{thm:Discrep} implies that
$$
\max_{1 \le a \le p-1} \left|\sum_{n=H+1}^{H+N}
\ep\( a\prod_{i=1}^\ell  n_i!\) \right|  \ll  N^{-\ell/4+r/2 - 1 +2^{-r}} 
p^{(\ell -2r+4)/8} (\log p)^{(\ell -2r+ 4)/4}
$$
for $\ell \ge 2r \ge 1$.

Let $G_{\ell}(a,N)$ be the number of solutions to the congruence:
$$
\prod_{i=1}^\ell n_i! \equiv a \pmod p,
$$
in positive integers $n_1, \ldots, n_\ell$ with
$$
     \sum_{i=1}^\ell n_i = N.
$$

It has been shown in~\cite{LS} that for any $\varepsilon$ and sufficiently
large $p$, $G_{\ell}(a,N) > 0$ provided that $\ell \ge p^\varepsilon$ and
$N -\ell >  p^{1/2+\varepsilon}$. In~\cite{GL}, the same result
has been obtained under a much weaker condition
     $N -\ell >  p^{1/4+\varepsilon}$. Here,
concentrate on the value of $\ell$ and show that
it can be taken as $\ell = O(1)$ provided $N  >  p^{1/2+\varepsilon}$.

\begin{theorem}
\label{thm:Fixed Sum Product}
For any fixed   integer  $\ell \ge 1$ and any
integer  $N$   with  $1 \le N < p/\ell, $
the following bound holds:
$$
\max_{a \in \F_p^*} \left|G_\ell(a,N) - \frac{1}{p-1} {N-1 \choose
\ell-1}\right|
\ll
     N^{3\ell/4} p^{(\ell + 6)/8} (\log p)^{(\ell - 2)/4}.
$$
\end{theorem}

\begin{proof}
For $a \not \equiv 0 \pmod p$, we have
$$
G_\ell(a,N)= \frac{1}{p-1}
     \sum_{\substack{n_1, \ldots, n_t  \ge 1\\ n_1 +\ldots + n_\ell = N}}\,
\sum_{\chi \in \cX} \chi\( a^{-1} \prod_{i=1}^\ell n_i!\),
$$
where the sum is taken over all  multiplicative characters $\chi$ modulo $p$.
Separating the contribution from the principal character $\chi_0$, we obtain
$$
\left| G_\ell(a,N) - \frac{1}{p-1} {N-1 \choose \ell-1} \right| \le
\frac{1}{p-1}  R,
$$
where
\begin{eqnarray*}
R & = & \sum_{\chi \in \cX^*}\chi(a^{-1})\sum_{\substack{n_1,
\ldots, n_\ell \ge 1\\ n_1
+\ldots + n_\ell = s}}
     \chi\( \prod_{i=1}^\ell n_i!\) \\
&=& \sum_{\chi \in \cX^*} \chi(a^{-1}) \sum_{n_1, \ldots, n_\ell =1}^\ell
     \chi\(  \prod_{i=1}^\ell n_i!\)    \frac{1}{p} \sum_{c=0}^{p-1}
\ep( c( n_1 +\ldots + n_\ell -N)),
\end{eqnarray*}
because if $\ell N < p$ then the congruence $n_1 +\ldots + n_\ell
\equiv s \pmod p$
with $1 \le n_1, \ldots,  n_\ell \le  N$
is equivalent to the equation $n_1 +\ldots + n_\ell = s$.
Therefore,
$$
R \le   \frac{1}{p} \sum_{c=0}^{p-1} \sum_{\chi \in \cX^*}
\left| \sum_{n=1}^\ell
     \chi\(n!\)  \ep( c n)\right|^\ell .
$$
Arguing as in the proof of Theorem~\ref{thm:Two Products},
we  derive  the result follows from~\eqref{eq:T and I} and 
Theorems~\ref{thm:Char-Indiv} and~\ref{thm:AverValues}.
\end{proof}

For example, for any fixed $\varepsilon> 0$ and $p/\ell > N \ge
p^{1/2 +\varepsilon}$,
we have
$$
G_\ell(a,N)  = \frac{1}{p-1} {N-1 \choose
\ell-1} (1 + o(1))
$$
for every fixed $\ell > \varepsilon^{-1} + 4$.

We remark that one can easily drop the condition $N <p/\ell$ in
Theorem~\ref{thm:Fixed Sum Product}.

Let  $F(a,H,N) =F_1 (a,H,N)$ be the number of solutions
of the congruence $n! \equiv  a \pmod p$,  $H +1 \le n \le H+N$.

\begin{theorem}
\label{thm:Dist of Values}
Let $H$ and $N$ be integers with  $0 \le H< H+ N < p$.
Then  following bound holds:
$$
\max_{a \in \F_p^*} F(a,H,N) \ll N^{2/3}.
$$
\end{theorem}

\begin{proof} Let $K>0$ be a parameter to be chosen later. Let
$$\cA=
\{H+ 1\le  n \le H+N \ |\  a\equiv n!\pmod p\}= \cA_1\cup \cA_2,
$$
where
$$
\cA_1=\{n \in \cA \ | \ |n - m |>K \ \text{for all}\ m \ne n,
    m \in \cA\}\qquad \text{and} \qquad \cA_2=\cA\backslash \cA_1.
$$
It is clear that $\#\cA_1\ll N/K$. Assume now that $n\in \cA_2$.
Then there exists a nonzero integer $k$ with $|k|\le K$ and such that
$n!\equiv (n+k)!\pmod p$.
For each $k$, the above relation leads to a polynomial
congruence  in $n$ of degree $|k|$ and therefore it has at most
$|k| \le K$ solutions $n$. Summing up over all values of $k$ with
$|k|\le K$, we
get that $\#\cA_2\ll K^2$. Thus,
$$F(a,H,N) = \#\cA  =  \#\cA_1  +  \#\cA_2 \ll \frac{N} {K}+K^2,$$
and choosing $K=\fl{N^{1/3}}$ we get the desired inequality. \end{proof}

We have seen that the Wilson theorem immediately implies
the inequality $V_2(0, p-1) \ge (p-1)/2$. We now show
that this bound can be slightly improved.

\begin{theorem}
\label{thm:ValueSet-2}
The  following bound holds:
$$
V_2(0, p-1) \ge \frac{5}{8}p + O(p^{1/2} \log^2 p).
$$
\end{theorem}

\begin{proof}  By~\eqref{Wilson Repr}, we see that if
$a  \equiv b^{-1} \pmod p$,
   $1 \le b \le p-1$  is odd then $a \in V_2(0, p-1) $.

By ~\eqref{Wilson Thm}, we see that if
$a \equiv c ^{-1}(c+1)^{-1} \pmod p$ with some even $c = 2u$,  $1
\le c \le p-3$,
then $a \in  V_2(0, p-1) $ too.
Thus each such  $c$ which corresponds to an even $b = 2v$ in the above
representation, contributes one new element to $V_2(0, p-1)$.
It is also clear that no more than two distinct values of $c$ can
contribute the
same element.

Therefore,  $V_2(0, p-1) \ge (p-1)/2 + W/2$,
where  $W$ is the number of solutions of the congruence
$$
2u  (2u + 1) \equiv 2v  \pmod p, \qquad 0 \le u,v \le (p-3)/2.
$$
The Weil bound
 yields $ W = p/4 + O(p^{1/2} \log^2p)$
(see~\cite{Chalk}), which concludes the proof.
\end{proof}

We remark that Theorem~\ref{thm:ValueSet-2}     immediately 
implies that for every  integer $a$ there exists a
representation $ a\equiv n_1! n_2! + n_3! n_4! \pmod p$
with some positive integers $n_1, n_2, n_3, n_4$. 

\section{Concluding Remarks}

Most of our  results hold in more general settings. For example,
let $m\ge 1$ be any fixed positive integer and put
$$
T(m,\chi,f,H,N)= \sum_{n = H+1}^{H+N} \chi\(\prod_{\nu=1}^m (n+\nu-1)!\)
\ep(f(n)).
$$
Then Theorem~\ref{thm:Char-Indiv} holds with $T(\chi,f,H,N)$ replaced
by $T(m,\chi,f,H,N)$. In particular, if we write $Q(m,H,N)$ for the number
of $n=H+1,\ldots,H+N$ such that $n!,\ldots,(n+m-1)!$ are all primitive roots
modulo $p$, then the estimate
$$
Q(m,H,N) = N  \(\frac{\varphi(p-1)}{p-1} \)^m
+ O\( N^{1-1/2}  p^{(\ell+2)/4\ell(\ell+1)} \)
$$
holds for any fixed integer $\ell\ge 1$.

Let  $\cQ$  be the set
of all  distinct prime
divisors of $p-1$. For a set $\cR \subseteq \cQ$,
we denote by $T(\cR,H,N)$ the number of $n =H+1, \ldots, H+N$
such that for every $q \in \cQ$, $n!$ is a $q$th power residue modulo $p$
if and only if $q\in \cR$. Then the estimate
$$
T(\cR,H,N) = N \prod_{q \in \cR} \frac{q-1}{q}
\prod_{q \in \cQ \backslash \cR} \frac{1}{q}
+ O\( N^{1 -1/2\ell}  p^{(\ell+2)/4\ell(\ell+1)} \)
$$
holds for any fixed integer $\ell\ge 1$.

Techniques of the present paper apply also
to the sequences
   $$
\binom{2n}{n} = \frac{(2n)!}{(n!)^2}, \qquad (2n+1)!! = 1\cdot 3
\ldots \cdot (2n+1),
$$
and many others, as well as their combinations.

Also, with some minor adjustments, our  methods can
be used   to obtain  similar, albeit
somewhat weaker results for composite moduli. In this setup, our
basic tools such as the Weil bound and the Lagrange theorem,
have to be replaced with their analogues in residue rings modulo
a composite number. See, for example,~\cite{CoZh} for bounds of character
sums, and~\cite{KoSte} for
bounds on the number of small solutions of polynomial
congruences.

While the results of the present paper represent some progress towards
better
understanding the behaviour of $n!$ modulo $p$, there are several
challenging
questions that deserve further investigation. For example, our
Theorem~\ref{thm:Dist of Values} gives a nontrivial upper bound on
$F(a,H,N)$, but we conjecture that this result is far from being
sharp. We do not have any nontrivial individual upper bounds for $S(a,H,N)$.

Certainly, studying $V_1(H,N)$ is of primal interest.
Trivially, we have $V_1(H,N) \ge (N-1)^{1/2}$
(to see this it is enough to recall that $n = n!/(n-1)!$),
but we have not been able to obtain any  better lower bound.
In the opposite direction,
answering a question of Erd\H os, Rokowska and Schinzel~\cite{RS}
have showed that if the residues of $2!, 3!, \ldots, (p-1)!$ modulo $p$
are all distinct, then the missing residue must be that of $-((p-1)/2)!$,
that $p\equiv 5\pmod
8$, and that no such $p$  exists in the interval $[7,1000]$, but it
does not seem to be even
known  that there can be only finitely many such $p$, or, equivalently,
that $V_1(0,p-1)=p-2$ can happen only for finitely many values of 
the prime $p$.

It is very tempting to try to generalize the proof of 
Theorem~\ref{thm:ValueSet-2}
and consider longer products $c(c+1)\ldots (c+m)$. This may lead to
  an improvement of the constant $5/8$ of Theorem~\ref{thm:ValueSet-2}.
However, to implement this strategy one has to study
in detail image sets of such polynomials (and their overlaps),
which may involve rather complicated machinery.

\end{document}